\documentclass[12pt,a4paper]{article}
 \usepackage{amsmath,amstext,amssymb,amscd}

\usepackage{graphicx}

\begin{document}

\title{Majority rule  on rhombus tilings and Condorcet super-domains \\
\thanks{Supported in part by grant RFFI
20-010-00569-A). (Danilov, Karzanov).}}

\author{Vladimir I. Danilov\thanks{CEMI RAS, Nahimovskii prospect, 47, 117418 Moscow; email: vdanilov43@mail.ru.} 
 \and
Alexander V. Karzanov\thanks{CEMI RAS, Nahimovskii prospect, 47, 117418 Moscow; email: akarzanov7@gmail.com} 
 \and
Gleb A. Koshevoy\thanks{IITP RAS, Bolshoi Karetny 19, 127051 Moscow; email:
koshevoyga@gmail.com. }}

\date{\today}
\maketitle

\begin{abstract}
In this paper we consider a Condorcet domain (CD) formed by a rhombus tiling as a
voting design and consider a problem of aggregation  of voting designs using majority
rule. A Condorcet super-domain is a collection of CDs obtained from rhombus tilings on
a zonogone Z(n; 2) with the property that if voting designs (ballots) belong to this
collection, then the simple majority rule does not yield cycles. A study  of Condorcet
super-domains and methods of constructing them form the main subject of this paper.
\end{abstract}

\section{Introduction}

A Condorcet domain is a set of linear orders on a finite set of candidates (alternatives in a voting) such that if all preferences of the  voters are linear orders belonging to this set, then the simple majority rule does not yield cycles. We use the abbreviation CD for a Condorcet domain. 

A CD is called \emph{maximal} if no new linear order can be added to this CD so that the extended set be a CD again.
A CD is called \emph{normal} if it contains the \emph{standard} linear order $\alpha = (1 < 2 <\cdots  < n)$ and the opposite order $\omega = (n < n - 1 <\cdots < 1)$, called the \emph{anti-natural} one.

In \cite{DK-13} it was shown that there is a bijection between the set of maximal normal CDs and
the set of maximal cliques in the Bruhat lattice on linear orders with respect to the so-called compatibility relation.
This relation on linear orders was introduced by Chameni-Nembua~\cite{CN}, it says that two linear orders are \emph{compatible} if both the intersection and union of   the sets of inversions in these orders determine linear orders as well. (Usually to define the Bruhat lattice and the compatibility relation, linear orders on $[n]$ are associated with their inversion sets, which are special subsets of the set $\Omega=\{(i,j), \, i<j, i,j\in [n]\}$.)
Moreover, each CD  is a distributive sublattice of the Bruhat lattice.

Abello~\cite{A}, and Galambos and Reiner~\cite{GR} constructed a maximal CD in  the form of the union of maximal chains in the Bruhat lattice. In~\cite{DKK-12} we proposed a construction of maximal CDs by use of rhombus tiling diagrams and showed that this construction unifies  the ones due to Abello, and Galambos and Reiner.

In this paper we consider a CD formed by a rhombus tiling as a {\em voting design} and consider a problem of aggregation voting designs using majority rule.  Thus, the problem of the non-acyclicity under the majority rule is shifted in a higher level, from linear orders to voting designs.

We define \emph{Condorcet super-domain} to be a collection of CDs obtained from rhombus tilings on a zonogone $Z(n;2)$ with the property that if voting designs (ballots) belong to this collection, then the simple majority rule does not yield cycles. We use the abbreviation CSD for a Condorcet super-domain. Accordingly, CSD is called \emph{maximal} if one cannot add to this CSD one more CD of the tiling type so that the extended collection by a CSD either.
A study of such CSDs and methods of constructing them form the main subject of this paper.

It turns out that constructions of maximal CSDs have similarities to the problem of constructing maximal CDs. Similar to linear orders, one can introduce a reasonable notion of inversions for rhombus tilings. As a result, we obtain a bijection between the set of rhombus tilings on the zonogone $Z(n;2)$ and the set of their inversions, viewed as special triples $(i,j,k)$ in $[n]$. The rhombus tilings form a poset generated by the inclusion relation on the inversion sets. This enables us to apply the majority (or the median) rule for aggregating rhombus tilings. In general, like in the case of linear orders, this aggregation may result in a set of inversions which corresponds to none of the tilings. Nevertheless, when dealing with CSDs of our interest the aggregating procedure always produces a tiling. 
 
We show in Theorem 2 that maximal chains in the poset of  inversion sets (ordered by inclusion) give rise to a maximal CSD. Another class of CSDs consists of symmetric CSDs characterised in Theorem 3. 
 
The paper is organised as follows. Section 2 gives an overview of basic results on tilings needed to us. In Section 3 we define CSD and show that the problem of aggregation of tilings of a given domain reduces to the problem of aggregation for triples of  tilings from this domain. This gives a characterization CSDs in terms of the existence of medians in the graph of rhombus tilings. In Section~4 we obtain results on median Condorcet super-domains, viewed as analogs of results by Slinko and Puppe~\cite{Pup-Sli} on ordinary Condorcet domains. Section 5 introduces normal CSDs and a notion of compatibility (similar to the 
Chameni-Nembua relation [8] for linear orders) and shows a bijection between normal CSDs and cliques w.r.t. the compatibility relation. In the final sections, two classes of normal CSD are constrcuted. The first one is obtained by using \emph{cubillages}, three-dimensional analogs of rhombus tilings. The second class is formed by symmetric CSDs. In reality there exist more general, "mixed" CSDs, but we do not come into a study of these domains and restrict ourselves with only one example.

The theory of aggregation of tilings and Condorcet supper-domains developed in this paper reminds, in many aspects, the theory of aggregation of linear orders and Condorcet domains developed earlier (see, e.g.,~\cite{DKK-12}). This is not a coincidence by chance; these are samples in dimensions 1 and 2 of more general phenomena arising in the theory of cubillages of any dimension (for a survey on cubillages, see~\cite{DKK-19}). However, these generalizations are beyond our paper. 

\section{Rhombus tilings: basic facts}

We fix an integer $n\ge 2$ and refer to elements of the set $[n] = \{1,\ldots, n\}$ as \emph{colors}. 

In the upper half-plane $\mathbb R\times \mathbb R_{>0}$, we fix
$n$ vectors $\xi_1 ,\ldots,\xi_ n$ ordered clockwise around the origin $(0,0)$.
It is convenient to assume that these vectors have the same
length. The Minkowski sum of $n$ segments $[0,\xi_ i]$, $i=1,\ldots,n$,
forms a \emph{zonogon}; we denote it by $Z(n;2)$. In other words,
$Z(n;2)$ is the set of points $\sum _i a_i \xi_ i$ over all $0\le a
_i\le 1$. It is a center-symmetric $2n$-gon with the \emph{bottom}
vertex $b=(0,0)$ and the \emph{top} vertex $t=\xi_1 +\cdots+\xi_ n$. A
{\em tile} is a rhombus congruent to the sum of two segments $[0,\xi_i ]$ and
$[0,\xi_ j]$ for distinct $i,j\in[n]$; depending on the context, we call it an $i$- or $j$- or $ij$-\emph{tile}.

\begin{figure}[htb]
\begin{center}
\includegraphics[scale=.2]{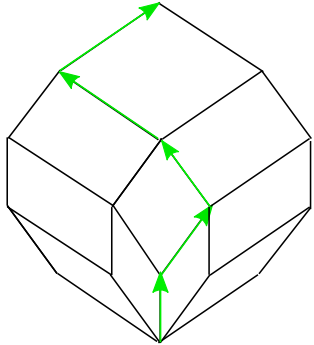}
\end{center}
 \caption{\small The standard tiling for $n=5$. The snake indicated by arrows gives the order $3<4<2<1<5$.} \label{fig:3}
 \end{figure}

A {\em rhombus tiling} (or simply a {\em tiling}) is a subdivision
$T$ of the zonogon $Z(n;2)$ into tiles which satisfy the
following condition: if two tiles intersect, then the
intersection consists of a common vertex or a common edge. Figure~\ref{fig:3} illustrates an example of tilings for $n= 5$.

Orienting the edges of $T$ upward, we obtain a 
planar directed graph (digraph) on the set of vertices of $T$, denoted as $G_T$. The tiles of $T$ are just the (inner two-dimensional) faces of $G_T$. An edge congruent to $\xi_i$ is called an edge of color $i$, or an $i$-\emph{edge}. 

For $i\in[n]$, the union of $i$-tiles in $T$ is called the $i$-{\em track}. One easily shows that
the $i$-tiles form a sequence in which any two consecutive tiles
have a common $i$-edge, and the first (last) tile contains the
$i$-edge lying on the left (resp. right) boundary of $Z(n;2)$.

By a {\em snake} of $T$, we mean a directed path in the digraph $G_T$ going from the bottom vertex $b$ to the top vertex $t$. 

Any snake contains exactly one $i$-edge, for
each $i$. So the sequence of ``colors'' of edges in a snake
constitutes a word $\sigma =i_1 \cdots i _n$, which is a linear
order on $[n]$. In what follows, we do not distinguish between
snakes $S$ and their corresponding linear orders $\sigma$,
denoting $S$ as $\mathcal S (\sigma )$ and  saying that the
linear order $\sigma $ is {\em compatible} with the tiling $T$.
The set of linear orders compatible with $T$ is denoted by $\Sigma
(T)$.

An important fact (see \cite{DKK-12}) is that the set $\Sigma(T)$ is a (maximal) Condorcet domain. So $\Sigma(T)$ can be considered as a good voting design, and the task of aggregation of such designs (tilings) looks reasonable (see the next section).

The set of tilings of the zonogon $Z (n, 2)$ is denoted by ${\bf T}_n$, or simply
${\bf T}$.  This set is endowed with several structures, resembling the set ${\bf L}={\bf L}_n$ of  linear orders on the set $[n]$.

\begin{figure}[htb]
\begin{center}
\includegraphics[scale=.3]{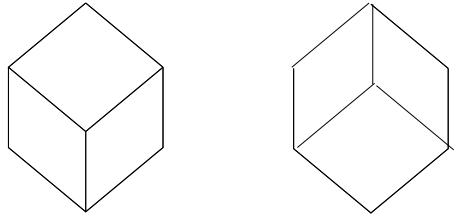}
\end{center}
 \caption{\small The standard tiling (left) and the anti-standard tiling (right) for $n=3$.}
 \label{fig:forn3}
 \end{figure}

1) For $n = 3$, there are only two tilings, see Fig.~\ref{fig:forn3}.

\noindent Here the left tiling is called \emph{standard}, and the right tiling  \emph{anti-standard}. 
For an arbitrary $n$, a tiling is called \emph{standard} if its reduction to any triple of colors is standard; the anti-standard tiling is defined symmetrically.

For a tiling $T$, one can form the antipodal tiling $T^\circ$, by doing the central symmetry transformation of  $Z(n; 2)$ and $T$ (relative to the center of the zonogon). Then the standard and anti-standard tilings are antipodal.
\medskip

2) In what follows, we will use an operation of reducing a tiling $T$ 
by a color $i\in [n]$. Consider the $i$-track $Q$ of $T$.  If we remove the interior of $Q$ and merge, in a natural way, the parts of $T$ lying above and below $Q$, we obtain the zonogon $Z(n-1;2)$ and a tiling of this zonogon, denoted as $T_{[n]-i}$ and called the \emph{reduction} of $T$ by the color $i$. Iterating this operation with other colors, we obtain the reduction of $T$ w.r.t. the collection of these colors.  For details, see, e.g. \cite{DKK-10}.

\begin{figure}[htb]
\begin{center}
\includegraphics[scale=.7]{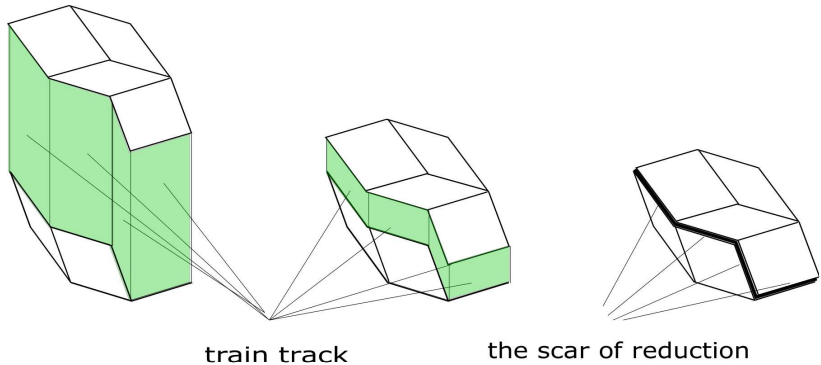}
\end{center}
\caption{\small A track (with color 4) and the corresponding reduction}
\label{fig:2}
\end{figure}

If $K$ is a subset in $[n]$, then for a tiling $T$, one can form the tiling of the zonogon $Z(|K|; 2)$, denoted as $T|_K$ or $T(K)$, by removing all colors not in $K$, i.e., $T|_K$ is the reduction of $T$ by $[n]-K$. This gives the reduction map
$$
T\mapsto T|_K.
$$
Clearly the standard tiling is reduced to the standard one, and similarly the anti-standard tiling.\medskip

3) \emph{Binary representation} of tilings. Denote by $\Lambda =\Lambda^3[n] =\tbinom{[n]}{3}$ the set of  triples $ijk$ of elements $i < j < k$ in $[n]$. For a tiling $T$, define the set $Inv(T)$ (of \emph{inversions}) consisting of the triples $ijk\in \Lambda$ for which the restriction $T|_{\{i,j,k\}}$ is anti-standard. In particular, $Inv(T_{st})=\emptyset$, and $Inv(T_{as})=\Lambda$. The mapping 
$$
Inv : {\bf T} \to 2^\Lambda
$$
is injective (in other words, inversions determine tilings; this is a relatively simple known fact), which allows us to speak of the \emph{binary representation}, or the \emph{binary form}, of a tiling.

It can be seen that for a subset $K$ of $[n]$, the following diagram is commutative:
$$
\begin{CD}
{\bf T} @>Inv>> 2^\Lambda \\
@VVV  @VVV \\
{\bf T}(K) @>Inv>> 2^{\Lambda (K)},
\end{CD}                                                  $$
where the first vertical arrow is the restriction map $|_K$ and the second one sends a subset $P$ of $[n]$ to $P\cap \Lambda^3(K)$. 

When we deal with an \emph{arbitrary} subset of $\Lambda$, we interpret it as a \emph{pseudo-tiling}. 

In order to decide if a pseudo-tiling is a tiling, that is, if a given subset of $\Lambda$ takes the form $Inv(T)$ for some $T\in {\bf T})$, one can use Ziegler's theorem; we review this later. 

It is clear that $Inv(T^\circ)=\Lambda  - Inv(T)$. Admitting some freedom, we will write $T\subset T'$ if $Inv(T)\subset Inv(T')$. The size $|Inv(T)|$ of  $Inv(T)$ is called the \emph{rank} of $T$. The standard tiling has rank 0, while the anti-standard one has rank $\tbinom{n}{3}=n(n-1)(n-2)/6$.\medskip

{\bf Remark}. There is another kind of binary representation of tilings using compatible orderings on $\Lambda^2([n])$, which goes back to Manin and Schechtman \cite{MSch}. Specifically, for any $i<j<k$, such an ordering
ranges the pairs in the triple ($ij$, $ik$, $jk$) as   $ij<ik<jk$ or $ij >ik>jk$.
 \medskip

4) \emph{Construction of the graph on $\mathbf T$.} Let us connect tilings $T,T'$ by a directed edge from $T$ to $T'$ if $Inv(T')$ is obtained from $Inv(T)$ by adding one element (one triple $ijk$).
It terms of tilings, this corresponds to a \emph{raising flip} in the hexagon with colors $ijk$, which replaces the standard tiling occurring in this hexagon by the anti-standard one.
Since an edge of the graph increases the rank of the tiling by one, the digraph $({\bf T}, \to)$ constructed in this way is 
acyclic, and therefore, the reflexive-transitive closure of the relation $\to$ determines a partial order $\preceq $ on {\bf T}. An important fact (\cite{FW}) is that this order $\preceq $ coincides with the order $\subset$ introduced for tilings (via their inversion sets) above. We should warn that a  pseudo-tiling that is the union or intersection of two tilings need not be a tiling. Moreover, the poset $({\bf T},\preceq)$ is not a lattice.

However, the structure of the underlying undirected graph for $({\bf T}, \to)$ is of importance to us. We denote this graph by $({\bf T}, \approx)$.
In case $n = 3$, $({\bf T}, \approx)$ is just a pair of vertices connected by an edge. For $n = 4$, $({\bf T}, \approx)$ is an octagon, a cycle with eight vertices. For $n = 5$, it is a graph with $62$ vertices, which is the skeleton graph of a three-dimensional zonotope (drawn in~\cite{FZ}, Figs. 2,3).

\begin{figure}[htb]
\begin{center}
\includegraphics[scale=.25]{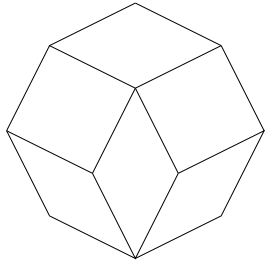}
\end{center}
 \caption{\small The standard tiling for $n=4$} 
 \label{fig:n=4}
 \end{figure}

Let us illustrate $({\bf T}, \approx)$ for $n = 4$. Standard tiling is viewed as in Fig.~\ref{fig:n=4}.

\begin{figure}[htb]
\begin{center}
\includegraphics[scale=0.7]{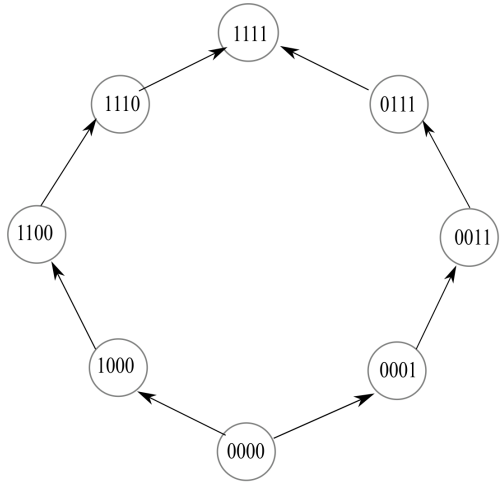}
\end{center}
 \caption{\small The octahedral  wheel} 
  \label{fig:8oct}
  \end{figure}

The other 7 tilings of the zonogon $Z(4, 2)$ are easy to draw as well. All tilings are naturally placed at the vertices of an octagon. The set $\Lambda^3[4]$ is the quadruple of triples $\{123; 124; 134; 234\}$,  and we accordingly encode the subsets of $\Lambda^3[4]$ by the corresponding strings (quadruples) consisting of zeros and ones. This encodes $2^4=16$ pseudo-tilings, and among these our $8$ tilings are as illustrated in Fig.~\ref{fig:8oct}. 

We observe that these 8 tilings are (encoded by) intervals in the set $[4]$, pressed to the left or to the right.
One can say that both ones and zeros should form a connected subset (interval) in the 4-string. (For this reason, they are also called \emph{bi-intervals}.) The central symmetry of the "octahedral  wheel" reflects the fact that for each tiling $T$, there is the ``opposite" tiling $T^\circ$. 
\medskip

\textbf{Theorem (Ziegler's criterion)}~\cite{Z}.  \emph{A pseudo-tiling $P$ is a tiling if and only if for any quadruple $F = ijkl$ of colors with $i < j < k < l$, the restriction $P|_F: = P \cap \Lambda^3(F)$ is a tiling of the zonogon $Z(F, 2)$}.
 \medskip
 
\noindent For a proof see also~\cite[Th.~15.8]{DKK-19}.

This criterion is convenient to check whether a pseudo-tiling is a tiling, and also to prove the equality of two tilings: it suffices to compare the quadruples. 

From now on we prefer to leave the geometric language and understand tilings as pseudo-tilings that satisfy Ziegler's criterion.
\medskip

\textbf{Example 1.} Let $I = [a\ldots b]$ be an interval in the chain $[n]$. Form the pseudo-tiling $T(I)$ consisting of all triples $ijk$ lying in $I$ (i.e., with $a \le i < k \le b$).  From Ziegler's criterion, it easily follows that $T(I)$ is a tiling of the zonogon $Z(n; 2)$, and the restriction $T(I)$ to the set of colors in $I$, $T(I)_I$, is the anti-standard tiling of the zonogon $Z(I, 2)$.
 \smallskip

One can also consider any tiling of the zonogon $Z(I; 2)$ as a subtiling of the ``large'' zonogon $Z(n; 2)$. The interested reader may try to construct geometrically the extension of the tiling with the color-set $I$ to the whole $[n]$. Such an extension is unique and gives the canonical embedding of $T(I)$ in $T$, which is ``inverse" to the reduction mapping $|_I$ . Note also that for any tiling $T\in {\bf T}$, the intersection $T\cap T(I)$ is again a tiling (and coincides with $T|_I$).

In the above constructions it is important that $I$ is an interval. Otherwise all related statements become incorrect in general.\medskip

\textbf{Example 2.} 
Recall that each tiling $T$ is associated with CD $\Sigma(T)$ consisting of linear orders on $[n]$ being 
%%. This domain was defined in geometric language as 
the set of snakes in the tiling $T$. Now we want to rephrase the definition  $\Sigma (T)$ in combinatorial language (especially since a similar translation will be used once again when constructing "cubilliage" super-domains). In other words, we want to give a ``purely combinatorial'' method to decide whether one or another linear order belongs to the domain $\Sigma (T)$.

Each triple $i<j<k$ of colors determines a triple of elements of the set $\Lambda^2[n]$ (formed by the pairs of elements of $[n]$), namely, the triple  $(ij,\, ik,\,jk)$. Note that these pairs are written in the lexicographic
order. For any linear order $<_\sigma$ on $[n]$, the set of its inversions $Inv(<_\sigma )$ intersects the triple  $(ij,\, ik,\,jk)$ either by a beginning interval or  by a terminal one. In the first case, we say that the order $<_\sigma$ is \emph{left} on the triple $ijk$, in the second case, that it is \emph{right}. In this terminology, \emph{the linear order $<_\sigma$ belongs to the domain $\Sigma (T)$ if and only if it is left in any triple $ijk$ of the tiling (regarded as a set of inversions), and it is right in any triple that doesn't belong to  the inversion set of $T$.}

Suppose, for example, that $ n = 3$ and the tiling $T$ is anti-standard. Its inversion set  consists of the single triple $123$. Therefore, the corresponding area $\Sigma(T)$ consists of four linear orders that are expressed via inversion sets as $\emptyset$, $\{12\}$, $\{12, 13\}$, and $\{12, 13, 23\}$, and in the usual form as $(1<2<3)$, $(1<3<2)$, $(3<1<2)$, and  $(3<2<1)$. These are exactly the snakes in the anti-standard tiling of the zonogon $Z(3, 2)$.

\section{Condorcet super-domains}

As is mentioned in Introduction, a CD arising from a tiling on $Z(n;2)$  can be interpreted as a voting design
(of a special type) on a set of alternatives on $[n]$. We have interested in  a problem of aggregation designs, specifically, in aggregation of CDs of  tiling type. 

Let  $V$ be a group of voters, and each voter $v\in V$ offers some tiling $T_v$ as a design. The problem is to choose from the set $\{T_v$, $v\in V\}$, the most suitable design for the group, or more precisely, to form a "compromising tiling", which would meet the wishes of the majority of voters.

The binary representation of tilings (see item 2 in the previous section) provides us with the most convenient representation of tiling domains to define the majority rule on voting designs.
Namely, for a set $\mathcal T(V):=\{T_v$, $v\in V\}$, we consider the following set of triples 
$$
\mathcal I(V)=\{ijk \in \Lambda, \textrm{the number of } v\in V \textrm{ with } ijk\in Inv(T_v)  \textrm{ is more than the half of } V\}.
$$

In what follows, we assume for simplicity that the number of voters is odd. Let  $sm (\mathcal T(V))$ denote the pseudo-tiling with the inversion set  $\mathcal I(V)$. 
This pseudo-tiling need not be a tiling (and thus it may not provide us with a voting design). 

For example, for $n = 3$, let us take the three tilings represented by $(1110)$, $(0000)$, and $(0111)$ (see Fig.~\ref{fig:8oct}). The simple majority rule $sm$  gives a pseudo-tiling $(0110)$ that is not a tiling. This phenomenon is of a similar flavor as the "Condorcet paradox" for the majority rule on linear orders.

Next we are interested in ``restricted'' domains for which such "paradoxes" do not appear. To differ from the case of linear orders, we call subsets of ${\bf T}$ \emph{super-domains}  (whereas the term ``domain'' is used for subsets of  the set ${\bf L}$ of linear orders, as before). 
 \medskip
 
\textbf{Definition.} A super-domain ${\bf D} \subset {\bf T}$ is  a \emph{Condorcet super-domain } (CSD) if for any odd set $V$ of "voters", the majority rule $sm$ applied to tilings $T_v$ in ${\bf D}$, $v\in V$, produces a tiling as well, that is, subset of ${\bf T}$.
 \medskip

\textbf{Proposition 1.} \emph{{\bf D} is a Condorcet super-domain if and only if for any quadruple $F = \{i < j < k < l\}$, the restriction of ${\bf D}$ to $F$ is a Condorcet super-domain.}\medskip

\emph{Proof}. The restriction of ${\bf D}$ to F consists of all tilings $T|_F$ where $T$ runs over ${\bf D}$. Direction ``only if'' in this statement is obvious because
$$
sm ((\mathcal T(V))|_F = sm ((T_v|_F, v\in V )).
$$
Direction ``if'' follows from the same formula and Ziegler's criterion, see the item 5 of Section~2. \hfill $\Box$

Thus, a key role in understanding and construing Condorcet super-domains plays the case $n = 4$, which we will analyze in more details. For obvious reasons, we may restrict ourselves by a characterization of maximal (by inclusion) CSD, since an arbitrary CSD is obtained as a subdomain of a maximal CSD. The next description gives all maximal CSDs for $n = 4$.  There are two classes of such CSDs, one consisting of 5-element sets, and the other consisting of 4-element sets.\medskip

I. \emph{A CSD of size 5 consists of five consecutive tilings on the octahedral wheel.} \medskip

II. \emph{A CSD of size 4 consists of two pairs of opposite tilings.}\medskip

\emph{In addition, for any collection constituted of tilings  of  one of these classes, $T_v$, $v\in V$, each of $T_v$ is of class I or of class II, the aggregate tiling $sm (\mathcal T(V))$ coincides with a tiling $T_v$ for some voter $v\in V$. }\medskip

A proof that the super-domains of these classes are CSDs is straightforward. For the opposite direction, we rely on one simple fact: if three tilings on the ``octahedral wheel'' are located at distances $2$, $3$ and $3$ from each other, then the aggregation of this triple is not a tiling. For example: for the tilings $(1100)$, $(0111)$ and $(0001)$ (regarded as inversion sets), the sum is $(1211)$; so the median is $(0100)$, which is not a tiling. The same phenomenon takes place for the other triples of tilings with the distances as indicated.
\hfill $\Box$
\medskip

From this description we get the following important result.\medskip
 
\textbf{Theorem 1.} \emph{A super-domain ${\bf D}$ is a CSD if and only if for any triple of tilings $T_1$, $T_2$, and $T_3$, their aggregation $sm (T_1; T_2; T_3)$ is a tiling.}\medskip

\emph{Proof}. In view of Proposition 1, it suffices to verify the assertion for $n = 4$. In the latter case, the assertion follows from the explicit description of CSDs given above. \hfill $\Box$\medskip

\textbf{Corollary.} \emph{ Let $\mathbf D$ be a CSD, $V$ be a set of  odd cardinality, and let $T \in sm({\bf D}^V)$. Then the super-domain ${\bf D}\cup \{T\}$  is a CSD as well. }\medskip

To prove this, we again may assume that $n = 4$. Then the result follows from the fact that for $n = 4$ , there is a voter $v\in V$ such that $sm (\mathcal T(V)) = T_v$.
\hfill $\Box$
 \medskip

In particular, if a CSD ${\bf D}$ is maximal, then the majority rule $sm$ sends ${\mathbf D} ^V$ to ${\mathbf D }$. We call a CSD $\mathbf D$ \emph{closed} if  $sm ({\bf D}^V) \subset {\bf D}$ 
(this notion is in the same vein as in \cite{Pup-Sli}).\medskip

\textbf{Remark.} We can consider a more general case. Namely, for a {\em majority system}, we can define the corresponding generalized majority rule. 
More precisely, let $V$ be a set of voters. A \emph{majority system} is meant to be a nonempty subset $\mathcal F\subset 2^V$ (whose elements are called ``big coalitions'') which satisfies the following two requirements:
 \smallskip

(1) \emph{monotonicity}: a superset of a ``big coalition'' is ``big coalition'' as well, and
\smallskip

(2) \emph{deciding}: a coalition is ``big'' if and only if its complement is not ``big''.
\medskip

Here are two examples of majority systems. The first one is a dictatorial system. The second one is a usual simple majority system (with an odd number of voters). Both examples belong to the class of weighted majority systems; note that there are majority systems that are expressed via weights of voters.
\smallskip

Given a majority system $\mathcal F$, one can define the rule to aggregate tilings $\{T_v, \  v \in V \}$ by:
$$
T_\mathcal F= \{ijk \in \Lambda \colon \textrm{ the set of } v\in V \textrm{ with } ijk\in Inv(T_v) \textrm{ is a ``big coalition'' in } V\}.
$$

In general, the quasi-tiling $T_\mathcal F$ need  not to be a tiling. However, if  $T_v\in \mathbf D$ for all $v\in V$ and some CSD {\bf D}, then the quasi-tiling $T_\mathcal F$ is indeed a tiling. This can be seen from Post-Monjardet's theorem \cite{Monj-78} on a representation of majority systems by means of medians, or by considering the case $n = 4$ in Theorem 1.

\section{Median super-domains and median graphs}

Theorem 1 says that in order to prove that a super-domain {\bf D} 
is CSD,  it suffices to verify the majority rule for any triple of tilings from {\bf D}. Note that the output of the majority rule applied to a triples can be thought of as a median. Therefore, it is worth to reformulate Theorem 1 in terms of medians.

More precisely, for pseudo-tilings $P,Q$ (regarded as subsets of $\Lambda$), we say that a pseudo-tiling $R$ \emph{lies between} $P$ and $Q$ if $P\cap Q \subset R \subset P\cup Q$. A \emph{median} of three pseudo-tilings $P_1,P_2,P_3$ is a pseudo-tiling $P$ such that $P$ is simultaneously between 
$P_1$ and $P_2$, between $P_2$ and $P_3$, and between $P_3$ and $P_1$. Then such a $P$, if exists, is exactly the set       
  $$
                          (P_1\cap P_2)\cup (P_2\cap
      P_3)\cup (P_3\cap P_1).
      $$
If $T_1,T_2,T_3$ are three tilings, then their median $T$, if exists, is a pseudo-tiling which need not be a tiling. %However, if $T$ is a tiling,  then $T$ is unique, denoted as $t-med(T_1,T_2,T_3)$. 

A super-domain $\bf D \subset \bf T$ is called a \emph{median} one if for any three tilings from {\bf D}, their median belongs to ${\mathbf D}$. In these terms, Theorem 1 can be reformulated as follows: 
  \smallskip
  
($\ast$) \emph{Each closed CSD is a median super-domain, and vice versa}.
 \smallskip

For tilings $T$, $T'$, the set of tilings which are between $T$ and $T'$ is denoted by $[T, T']$ and called the \emph{interval} of $T$ and $T'$. A subset (super-domain) ${\bf D}\subseteq {\bf T}$  is called \emph{convex} if for any $T, T'\in{\bf D}$, the interval $[T, T']$ lies in {\bf D}. It is clear that a super-domain {\bf D} is a median domain if and only if it possesses the (classical) Helly property for convex subsets. That is if a family of sets has a nonempty intersection for every triple of sets, then the whole family has a nonempty intersection.

A non-triviality occurs if one wishes to represent intervals geometrically or graphically. The set ${\bf T}$ can be endowed with a structure of graph (see item 4 in Section 4). Therefore, we can consider graph intervals and use the language of graph geodesics. Let us say that for  tilings $T,Q$, a tiling $R$ lies in a {\em geodesic} between $T$ and $Q$ if $R$ belongs to some shortest path from $T$ to $Q$ in the graph $\mathbf T$. We denote by $[T;Q]^g$ the geodesic interval between $T$ and $Q$. An interesting question is how this geodesic interval is  related to the the interval $[T,Q]$ defined earlier. 

One can ask this question a little differently. If the tilings $T$ and $T'$ are adjacent in the graph {\bf T} (i.e., they differ by a single triple $ijk$), then the interval $[T, T']$ consists exactly of $T$ and $T'$. Is the converse true? In other words, is it true that if there are no other tilings between $T$ and $T'$, then $T$ and $T'$ differ by one triple? (Equivalently, in light of Felsner-Weil' theorem \cite{FW}: is it true that $T$ and $T'$ are comparable by inclusion.)
\medskip

\textbf{Example 3.} Here we demonstrate a nontrivial example of CSDs for $n = 5$ (the case of $n = 4$ was discussed in detail above). This super-domain consists of 16 tilings which are placed at the vertices of the directed graph drawn in the picture below; the edges (arrows) correspond to raising flips, and triples $ijk$ related to flips are indicated on arrows  (where ``parallel'' arrows have the same triples). The black bottommost vertex corresponds to the standard tiling. To obtain a tiling at a vertex $v$ of the graph, one should collect the triples written on the arrows of a shortest path from the bottommost vertex to $v$.

\

\unitlength=.8mm
\special{em:linewidth 0.4pt}
\linethickness{0.4pt}
\begin{picture}(121.00,101.00)(-15,0)
\put(70.00,5.00){\circle*{2.00}}
\put(70.00,25.00){\circle{2.00}}
\put(85.00,40.00){\circle{2.00}}
\put(55.00,40.00){\circle{2.00}}
\put(70.00,55.00){\circle{2.00}}
\put(90.00,65.00){\circle{2.00}}
\put(50.00,65.00){\circle{2.00}}
\put(35.00,50.00){\circle{2.00}}
\put(105.00,50.00){\circle{2.00}}
\put(70.00,75.00){\circle{2.00}}
\put(90.00,85.00){\circle{2.00}}
\put(50.00,85.00){\circle{2.00}}
\put(100.00,100.00){\circle{2.00}}
\put(40.00,100.00){\circle{2.00}}
\put(20.00,45.00){\circle{2.00}}
\put(120.00,45.00){\circle{2.00}}
\put(70.00,6.00){\vector(0,1){18.00}}
\put(71.00,26.00){\vector(1,1){13.00}}
\put(69.00,26.00){\vector(-1,1){13.00}}
\put(56.00,41.00){\vector(1,1){13.00}}
\put(84.00,41.00){\vector(-1,1){13.00}}
\put(54.00,41.00){\vector(-2,1){18.00}}
\put(86.00,41.00){\vector(2,1){18.00}}
\put(104.00,51.00){\vector(-1,1){13.00}}
\put(71.00,56.00){\vector(2,1){18.00}}
\put(69.00,56.00){\vector(-2,1){18.00}}
\put(36.00,51.00){\vector(1,1){13.00}}
\put(70.00,56.00){\vector(0,1){18.00}}
\put(50.00,66.00){\vector(0,1){18.00}}
\put(90.00,66.00){\vector(0,1){18.00}}
\put(71.00,76.00){\vector(2,1){18.00}}
\put(69.00,75.00){\vector(-2,1){18.00}}
\put(49.00,86.00){\vector(-2,3){9}}
\put(90.00,86.00){\vector(2,3){9.00}}
\put(106.00,50.00){\vector(3,-1){13.00}}
\put(34.00,50.00){\vector(-3,-1){13.00}}
\put(74.00,15.00){\makebox(0,0)[cc]{234}}
\put(81.00,32.00){\makebox(0,0)[cc]{235}}
\put(57.00,32.00){\makebox(0,0)[cc]{134}}
\put(99.00,44.00){\makebox(0,0)[cc]{245}}
\put(39.00,44.00){\makebox(0,0)[cc]{124}}
\put(74.00,66.00){\makebox(0,0)[cc]{135}}
\put(98.00,92.00){\makebox(0,0)[cc]{145}}
\put(41.00,92.00){\makebox(0,0)[cc]{125}}
\put(30.00,45.00){\makebox(0,0)[cc]{123}}
\put(112.00,45.00){\makebox(0,0)[cc]{345}}
\end{picture}

One can directly check that this set of 16 tilings forms a CSD; we guess that this CSD is maximal. Also one can check that the underlying undirected graph is a median graph (recall that a graph is called a median one if every triple of its vertices has a unique median).  For example, consider the following three tilings (regarded as inversion sets):
$T_1 = \{234; 235; 245\}$, 
$T_2 = \{234; 134; 124\}$, and
$T_3 = \{234; 235; 134; 135\}$.

The median of these sets is $\{234; 134; 235\}$ and the corresponding  tiling is the central vertex of the graph.
Note that this CSD contains the standard tiling, but does not contain the anti-standard one. The latter cannot be added. For otherwise the median  of $T_1$, $T_2$, and the anti-standard tiling is equal to $T_1\cup T_2 = \{234, 235, 245, 134, 124\}$, but this quasi-tiling  is not a tiling. To see this, consider the {\ em stick} $\{124, 125, 145, 245\}$ for the quadruple $1245$. Our set intersects this stick by two triples $124$ and $245$, which stand on the edges of the stick, and therefore does not  satisfy Ziegler's criterion.

%------------------ Sect. 5
\section{Normal super-domains}

In general, it is difficult to say anything explicit about the structure of maximum CSDs. For this reason, we restrict ourselves by considering the so-called normal CSDs.\medskip

\textbf{Definition.} A super-domain {\bf D} is called \emph{normal} (borrowing the term used for domains in {\bf L}) if it contains both the standard and anti-standard tilings. \medskip

For such super-domains, it is convenient to use the following  compatibility relation (analogous to the compatibility relation due to  Chameni-Nembua for linear orders).\medskip

\textbf{Definition.} Tilings  $T$ and $T'$ are called \emph{compatible} if both $T\cap T'$ and $T\cup T'$ are tilings, denoting this relation by $\sim$.\medskip

For example, if tilings $T,T'$ are comparable (by inclusion), then they are compatible. In particular, any tiling is compatible with the standard tiling and the anti-standard tiling. It is easy to see that  any tiling $T$ is comparable and therefore compatible with its opposite tiling $T^\circ$. If tilings $T,T'$ are compatible, then $T^\circ$ and $T'^\circ$ are compatible as well.

Let us explain the relation $\sim$ using the octagonal wheel ${\bf T}_4$ ?? drawn in Fig.~5. All tilings on the right branch of this wheel are compatible to each other, as well as all tilings on the left branch. As to the pairs of tilings in different branches, only the opposite tilings are compatible. This description easily follows from the following \medskip

\textbf{Lemma 1.} \emph{If tilings $T, Q, R$ are pairwise compatible, then $R$ is compatible with $T\cap Q$ and $T\cup Q$, and the equality $R \cap (T \cup Q) = (R \cap T) \cup (R \cap Q)$ holds true.}\medskip

{\em Proof}.
It  suffices to check this for $n = 4$, i.e., for the octagonal wheel. If $T\subset Q$, say, then the statements are trivial. If $T$ and $Q$ are opposite, then $T \cup Q$ is the anti-standard tiling, and $T \cap  Q$ is the standard one, and the statements follow.\hfill $\Box$
\medskip

A super-domain {\bf D} is called a {\em clique} if all elements of $\mathbf D$ are pairwise compatible.\medskip

\textbf{Corollary.} \emph{Let {\bf D} be a maximal clique in the graph $(\bf T, \sim)$. Then ${\bf D}$  is a distributive sublattice in the poset $(\bf T, \subset)$.}\medskip

\textbf{Proposition 2.} \emph{Let {\bf D} be a normal super-domain. Then the following assertions are equivalent:}

(i) {\bf D} \emph{is a clique}; 

(ii) {\bf D} \emph{is a CSD}.\medskip

{\em Proof}. It suffices to consider the case  $n = 4$, and in this case 
the equivalence of (i) and (ii) is straightforward. 
\hfill $\Box$
\medskip

In particular, the maximal normal CSDs are exactly the maximal cliques in $({\bf T}, \sim)$. By the way, one can  aggregate tilings in $\mathcal T(V)=(T_v, \  v \in V)$ by the explicit formula
                                   $$
                            sm(\mathcal T(V))=\cup _M(\cap _{v\in      M}T_v),
                                     $$
where $M$ runs over the collection of subsets in $V$ of size greater than $|V|/2$ ("big coalitions"). The same formula works in the case of a generalized majority.

Using Proposition 2 and Corollary, one can construct various normal CSD. 
Suppose we are given a set ${\bf D}$ of tilings which satisfies the following properties:

a) {\bf D} is closed with respect to the operations $\cap$ and $\cup$;

b) ${\bf D}$ contains both $\emptyset$ and $\Lambda$.

Then ${\bf D}$ is a normal Condorcet super-domain.

In what follows we discuss two large classes of normal CSDs that are construct based on this principle. For the first class, we take chains in the graph ${\bf T}$; for the second one, we consider Boolean sublattices.

\section{Cubillage CSDs}

Any chain of tilings that begins at the standard (empty) tiling and ends at the anti-standard one (viz. $\Lambda$) gives a normal Condorcet super-domain. To get a larger super-domain, we need to take a ``non-extendable'' (maximal) chain of tilings $T_0 = \emptyset,T_1,\ldots,T_k$ with $k=\binom{n}{3}$, that is, a chain consisting of $\binom{n}{3}+1$ different tilings. Each subsequent tiling $T_{i+1}$ is obtained from the previous one $T_i$ by adding one triple from $\Lambda$. In other words, a maximal chain of tilings defines (and is defined by) a sequence $(t_0,t_1, \ldots , t_k)$ of different elements of $\Lambda$ forming a linear order on $\Lambda$. Clearly this linear order is not arbitrary; at each step a tiling has to be duly formed. In other words, at each step one should obey Ziegler's criterion. The following definition formulates this in terms of orders on $\Lambda$.\medskip

\textbf{Definition} \cite{MSch}. The order (partial or linear) on $\Lambda= \Lambda^3 [n]$ is called \emph{admissible} if its restriction to any quadruple $\Lambda^3(\{i; j; k; l\})$ is lexicographic or anti-lexicographic.\medskip

Ziegler's criterion asserts that if $T$ is a \emph{stack} (ideal) of an admissible order on $\Lambda$, then $T$ is a tiling. In particular, an admissible linear order gives a non-extendable chain of tilings, thus forming a CSD. Note, however, that CSDs constructed in this way are, as a rule, not maximal. To get larger
Condorcet super-domains, we are forced to use partial admissible orders. In  geometric terms, this would lead us to appealing to the so-called cubillages of three-dimensional cyclic zonotopes. We, however, prefer to use the combinatorial language at this point (referring the curious reader to \cite{DKK-19} for a survey on cubillages).

\begin{figure}[htb]
\begin{center}
\includegraphics[scale=.34]{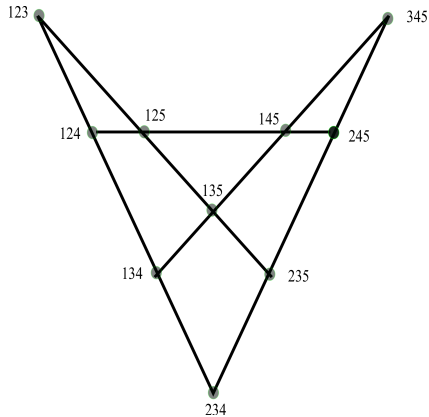}
\end{center}
 \caption{\small Grassmannian $\Lambda$ for $n=4$ covered by sticks. } \label{fig:palki}
 \end{figure}

For each quadruple $F = ijkl$ (where $i < j <k < l$), we form the sequence of four triples $ijk$, $ijl$, $ikl$, $jkl$, that we call a \emph{stick}. Then the set $\Lambda$ is covered by such sticks, and every triple (an element of $\Lambda$) is contained in $n-3$ sticks.  Figure~\ref{fig:palki} illustrates the case $n = 5$. (An illustration for $n = 6$ is given in \cite{Z}.)

Every stick has two distinguished linear orders: the lexicographic one: $ijk < ijl < ikl < jkl$, and the reverse to it, the anti-lexicographic one. In the lexicographic case, we say that the stick has the \emph{direct orientation}, and in the anti-lexicographic one, the \emph{reverse orientation}. Then an admissible order on set $\Lambda$ determines one or another  orientation on each stick. Conversely, orientations on the sticks determine an admissible order if  the obtained graph is acyclic (!); in fact, in this case we deal with a (partial) admissible order on $\Lambda$. Once having an admissible order $\preceq$ on $\Lambda$, we can form a super-domain ${\bf T} (\preceq)$ consisting of all stack-ideals (and hence tilings) with respect to this order $\preceq$. More precisely, a tiling $T$ is included in this super-domain if its restriction to (intersection with) every stick consists of an initial segment of the stick when it is oriented directly, or an end segment when it is oriented reversely. In other words, we require the consistency w.r.t. stacks only for the restrictions on every stick. Since the intersection and union of stacks are stacks again, all tilings in this domain are compatible to each other. Therefore (see Proposition 2), the super-domain ${\bf T} (\preceq)$ is a CSD and, moreover, it is a normal one. Thus, on this way we obtain a plenty of CSDs. They are maximal ones, as the following theorem shows.\medskip

{\bf Theorem 2.} \emph{For any admissible order $\preceq$ on $\Lambda$, the CSD ${\bf T}(\preceq)$ is maximal by inclusion.}\medskip

{\em Proof}. Assume this is not the case, and let $T$ be a tiling that does not belong to the CSD ${\bf T}(\preceq)$, but  ${\bf T}(\preceq)\cup \{T\}$ is a CSD. Since ${\bf T}(\preceq)$ is normal, this is equivalent to the property that $T$ is compatible with any tiling in ${\bf T}(\preceq)$.

Since the tiling $T$ cannot be added to ${\bf T}(\preceq)$, there is a stick (say, oriented directly) whose intersection with $T$ is not an initial segment of this stick. Then this intersection a nonempty end segment. The latter is a segment of the form (0001), (0011), or (0111). Since the restrictions of tilings from ${\bf T}(\preceq)$ give the whole set of initial segments in
this stick (which is seen by considering a maximal chain in ${\bf T}(\preceq)$), we can easily find a tiling in the cubillage system ${\bf T}(\preceq)$ that is incompatible with $T$.
\hfill $\Box$

\section{Symmetric CSDs}

Another extreme case concerns Boolean lattices which are formed by pairs of opposite to each other tilings $T$ and $T^\circ$. We call such CSDs  \emph{symmetric}. As before, we will consider merely normal CSDs. In addition, we  assume that each CSD we deal with is closed. Then the tilings of  such a symmetric CSD ${\bf D}$ form a distributive lattice of subsets of $\Lambda$, i.e., they are closed with respect to taking the intersection and union. Since the super-domain ${\bf D}$ is symmetric, it is closed with respect to taking the complement either. So it is a \emph{Boolean ring} of subsets of $\Lambda$. Let $T_1,\ldots,T_m$  be the minimal non-empty elements of ${\bf D}$; then they are pairwise disjoint and their union is the whole $\Lambda$. 

One can go in the reverse direction. Consider a partition of $\Lambda$ into $m$
(pairwise disjoint) tilings:
$$
            \Lambda      =T_1\sqcup \ldots\sqcup T_m
$$
(regarded as subsets of $\Lambda$). For a subset $S\subseteq [m]$, define $T_S = \cup_{s\in S} T_s$.\medskip

\textbf{Lemma 2.} \emph{For any $S\subseteq [m]$, $T_S$ is a tiling.}\medskip

{\em Proof}. As before, it suffices to consider the case $n = 4$. Then it is easy to see that all sets $T_s$ in the partition are empty, except for two sets, say, $T_1$ and $T_2$, which are opposite to each other. This implies that $T_S$ is either empty or the whole $\Lambda^3[4]$ or $T_1$ or $T_2$.
\hfill $\Box$
\medskip

In view of this, a super-domain ${\bf D}$ formed by all tilings $T_S$, where $S\subseteq[m]$, is a CSD and, moreover, it is a symmetric CSD (since the opposite to $T_S$ is $T_{[m]\setminus S}$). The size of this CSD is equal to $2^m$ (if all sets $T_s$ are nonempty). 

 We assert that $m\le n-2$. This follows from\medskip
 
\textbf{Lemma 3.} \emph{For any nontrivial (nonempty) tiling $T$, its inversion set  contains a dense triple.}\medskip

{\em Proof}. Define the \emph{amplitude} of a triple $ijk$, $i < j < k$,  to be the number $k - i$. Since $T$ is nonempty, it contains some triple $\tau=ijk$; assume that $\tau$ has the minimal amplitude among the triples in $T$. We assert that its amplitude is 2, i.e., $\tau$ is an interval. For otherwise, there is a `hole' $\ell$ inside the interval $[i..k]$; let for definiteness $j < \ell < k$. Consider the quadruple  $ijlk$ and the corresponding four triples (the stick) $ij\ell, ijk, i\ell k, j\ell k$. 
The second member $ijk$ of this stick belongs to $T$. Then at least one of $ij\ell$ or $j\ell k$ must belong to $T$ as well (by Ziegler's criterion). But each of these triples has a smaller amplitude, contrary to the assumption. This implies the lemma.\hfill $\Box$
\medskip

\textbf{Corollary.} $ m \le n-2$.\medskip

Thus, the size of any symmetric CSD does not exceed $2^{n-2}$. And if $m = n-2$, then we get a maximal CSD.  Moreover, the maximality here is achieved not only among symmetric CSDs but among the variety of all CSDs (see remarks on the compatibility definition). 

A natural problem is how to construct a partition of $\Lambda$ into tilings more explicitly? Especially partitions into $m = n-2$ tilings. This is discussed in the next section.

      \section{Construction of symmetric CSDs}
      
It was seen above that an important role in the analysis of tilings is played by``interval'' triples. A triple $ijk$ is called \emph{dense} if it forms an interval, i.e., $i=j-1$ and $k=j+1$. We also denote such a triple as $ext(j)$, where $j$ can take value from 2 through $n-1$. By Lemma 3, any nontrivial tiling $T$ contains a dense triple. In fact, from the proof of that lemma one can obtain a sharper assertion, as follows.\medskip

\textbf{Lemma 4.} \emph{If a tiling $T$ contains a triple $ijk$, then $T$ contains also a dense triple $ext(\ell)$ with $i < \ell < k$, i.e.,  a dense triple lying in the interval $[i\ldots k]$.}\medskip

\textbf{Definition.} The \emph{basis} of a tiling $T$ is the set $B (T)$ of dense triples belonging to $T$.\medskip

So $B$ may be thought of as a map of ${\bf T}$ to $[2\ldots n-1]$. It is far to be injective: one can construct many tilings with the same basis. An important  statement is that there exists a minimal tiling among those. Before proving this, it is useful to consider the case when the basis is an interval $I\subseteq [2\ldots n - 1]$).

For such an $I$, let $ext(I)$ denote corresponding  interval in $[n]$. Let $T(ext(I))$ be the tiling introduced in Example~1 from Section~2 (consisting of all triples $ijk$ that lie in $ext (I)$). Then the basis $B(T(ext (I)))$ of the tiling $T(ext (I))$ coincides with $I$. We assert that $T(ext(I))$ is the minimal tiling with basis $I$. \medskip

\textbf{Lemma 5.} \emph{Let $T$ be a tiling and let $B(T)$ contain  an interval $I$. Then the inversion set of $T(ext(I))$ is a subset of that of $T$.}\medskip

\emph{Proof}. Let $ijk$ be a triple of $T (ext(I))$. % (i.e., elements $i,j,k$ belong to $ext(I)$).?? 
We show that this triple lies in $T$. This is proved by induction
on the amplitude of  triples. If the amplitude is two, i.e., it is a dense triple, then $j\in I\subset B(T)$ and $ijk\in T$. So we may assume that the amplitude of $ijk$ is larger than two and that all triples in $T (ext(I))$ with smaller amplitudes belong to $T$. Since $ijk$ has the amplitude larger two, it has a ``hole'' $\ell$, say, $i < \ell < j < k$. Now form the corresponding stick of triples $i\ell j,\,i\ell k,\,ijk,\,\ell jk$. Both ends $i\ell j$ and $\ell jk$ of this stick  belong to $T(ext (I))$ (since they lie in $ext (I)$). Also they have smaller amplitudes. By induction, they belong to $T$. Then, according to Ziegler's criterion, the middle term $ijk$ also belongs to $T$. \hfill$\Box$\medskip

%%\textbf{Remark.} Although this is not needed to us in what follows, we can construct the minimum tiling $T_B$ (with the base $B$) for any $B\subset [2\ldots n-1]$. To do this, partition $B$ into a disjoint union of intervals $I_1, \ldots, I_m$ and form $T_B$ as the union of $T(ext(I_s))$ for $s = 1, \ldots, m$. It is easy to realize (using the same Ziegler's criterion) that this is a tiling with the basis $B$. Lemma~5 shows that it is a?? minimal tiling with this basis.
%%\medskip
%%Also it is easy to see that there is a?? maximal tiling $T^B$ with the basis $B$. Thus, the entire set (poset) $T$ is divided into $2^{n-2}$ intervals $[T_B; T^B]$.

Now we start to construct symmetric CSDs of the maximal size $2^{n-2}$. That is, we are going to partition $\Lambda$ into $n-2$ nonempty tilings $T_2, \ldots  ,T_{n-1}$. Take an arbitrary element $j$ in the interval $[2.. n-1]$). The constructions are slightly different when $j$ is 2 or $n-1$, and when $j$ is an inner element of the interval $[2\ldots n -1]$.
If $j = n-1$, say, we split $\Lambda$ into $T([1\ldots n-1])$ and its complement. And if  $j$ is inner, this $j$ gives the partition of the interval $[2.. n-1]$ into three subintervals: $\{j\}$, $I = [2.. j - 1]$ and $K = [j + 1.. n-1]$. Accordingly, $\Lambda$ is divided into the tilings $T(ext (I))$, $T (ext (K))$ and their complements. Here it is important to notice that the complement to $T(ext(I))\cup T(ext(K))$ is a tiling; we denote it as $T_j$. In fact, $T_j$ has inversions of the form of all possible triples $abc$ such that $a< j < b$. It is easy to check that Ziegler's criterion is satisfied for such a set of triples. Then one can continue this procedure with any element from $I$ or $K$, and so on. As a result, we get a partition $\Lambda = T ([2.. n-1])$ into $n-2$ tilings, as required.
 
Thus, everything is determined by the sequence $j_1,j_2, \ldots,  j_{n-2}$ of  different elements of $[2 \ldots n-1]$. For example, the picture below illustrates a partition for $n=7$ with the sequence $4, 5, 6, 2, 3$. It should be noted that different sequences can produce the same partition, yielding the same symmetric CSD. The matter is that if $j_1$ is an inner point of the interval $[2 \ldots n-1]$, then we do not care which of the two neighboring intervals ending at $j_1$ are chosen for the next step, and so on. 

\begin{center}
\includegraphics[scale=1]{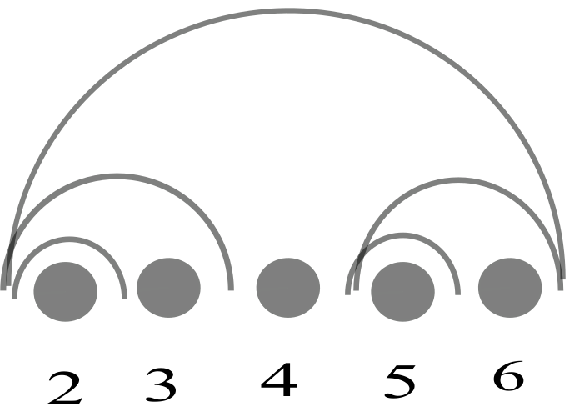}
 \end{center}

As a consequence, we obtain the following 
\medskip

\textbf{Theorem 3.} \emph{A linear order $j_1, \ldots,  j_{n-2}$ on the set $[2.. n-1]$ determines a symmetric CSD of the size $2^{n-2}$.} 
\medskip

%%We wished to say that this construction is universal and gives all partitions of $\Lambda$ into $n-2$ tilings. However, this is not true; here is an example for $n=5$. Namely, take $T_2 = \{123\}$, $T_3 = \{124,125,134,135, 234,235, 245\}$, $T_4 = \{145,345\}$. This is a partition of $\Lambda^3[5]$ into three tilings (which is checked straightforwardly).

      \end{document}